\documentclass[a4paper, 12pt]{article}
\usepackage{amssymb, amsmath, amsfonts}
\usepackage[english]{babel}
\usepackage{amsthm,amscd,amsfonts,latexsym,color,epsfig}
\usepackage{amssymb, epic, graphicx}
\begin{document}
\begin{center}

\sloppy

\textbf{\Large On boundary-value problems for a partial differential equation with Caputo and Bessel operators}

\bigskip

\bigskip

\centerline {  {\bf Praveen Agarwal}\footnote{Department of Mathematics, Anand International College of Engineering, Jaipur-303012, India\\ E-mail: goyal.praveen2011@gmail.com} {\bf Erkinjon Karimov}\footnote{ Institute of Mathematics, National University of Uzbekistan, Tashkent,100125, Uzbekistan
 \\ E-mail:erkinjon@gmail.com } {\bf Murat Mamchuev}\footnote{Department of Theoretical and Mathematical Physics, Institute of Applied Mathematics and Automation, Shortanova str. 89-A, Nalchik, 360000 Kabardino-Balkar Republic, Russia\\ E-mail: mamchuev@rambler.ru} {\bf and Michael Ruzhansky}\footnote{Department of Mathematics, Imperial College London, 180 Queen's Gate, London SW7 2AZ
United Kingdom\\
E-mail: m.ruzhansky@imperial.ac.uk}}

\end{center}

\bigskip

\begin{abstract}
In this work, we investigate a unique solvability of a direct and inverse source problem for a time-fractional partial differential equation with the Caputo and Bessel operators. Using spectral expansion method, we give explicit forms of solutions to formulated problems in terms of multinomial Mittag-Leffler and first kind Bessel functions.
\end{abstract}

\section{Introduction}

It is well-known that partial differential equations(PDEs) are playing a key role at constructing mathematical models for many real-life processes. Especially, during the last several decades, many applications of various kind of fractional differential equations became target of specialists due to both theoretic and practical reasons \cite{KST06}. Many kinds of boundary problems, including direct \cite{nakh} and inverse problems \cite{isak}, were formulated for different type of PDEs of integer order and with several fractional order differential operators.

Omitting survey part on fractional PDEs and inverse problems, we only note some works, closely related to the subject of this paper. In \cite{gorluch}, Luchko and Gorenflo studied fractional differential equation with Caputo fractional derivative and using operational method, solution of initial boundary problem for that equation was obtained in explicit form containing multinomial Mittag-Leffler function, some properties of which were obtained by Zhiyuan Li et al \cite{yamli}. Later, Yikan Liu \cite{liu} established strong maximum principle for fractional diffusion equations with multiple Caputo derivatives and investigated related inverse problem. We note as well the work by Daftardar-Gejji and Bhalikar \cite{gejbhal} where multi-term fractional diffusion-wave equation was considered and boundary-value problems for this equation were solved by the method of separation of variables.

There are many works on studying direct and inverse problems for time-fractional diffusion or diffusion-wave equations with the Caputo derivative, where operator acting on the space-variable is elliptic or uniformly elliptic. Depending on which operator is used in the space-variable, existence of classical or generalised solution is partly known, for instance, see works \cite{sakyam}, \cite{luch}, \cite{pengkar} and references therein. We also note works related to the Bessel operator. In \cite{masmes}, the initial inverse problem for heat equation with Bessel operator was investigated. In a recent work By Fatma Al-Musalhi et al \cite{fatma}, the authors studied inverse initial and inverse source problems for time-fractional diffusion equation with zero order Bessel operator. For review of different questions for equations with integer order derivatives we can refer to \cite{RS}.

A distinctive aspect of our problem is that with respect to the time-variable we have two Caputo derivatives and in a space variable our equation has a singularity at zero, hence it requires different regularity conditions from the given data.

First, we formulate direct and inverse-source problems and present preliminary information on Bessel equation, Fourier-Bessel series, and as well as on general solutions to the corresponding two-term fractional differential equation with Caputo derivative. At the end we formulate our result as two theorems.

\section{Formulation of problems}

We consider the equation
\begin{equation}\label{eq1}
  \partial_{0t}^{\alpha}u(x,t)+\lambda\, \partial_{0t}^{\beta}u(x,t)=u_{xx}(x,t)+\frac{1}{x}u_x(x,t)-\frac{\nu^2}{x^2}u(x,t)+f(x,t)
  \end{equation}
in a rectangular domain $\Omega=\left\{(x,t):\, 0<x<1,\,\,0<t<T\right\}$,
where $$0\leq \beta\leq \alpha,\,\,\beta<\alpha\le 1,$$ 
$\nu,\,\,T>0$. Here, in general,
\begin{equation*}
\partial_{0t}^{\alpha}g(t)=\left\{ \begin{aligned}
  & \frac{1}{\Gamma \left( k-\alpha  \right)}\int\limits_{0}^{t}{\frac{{g}^{(k)}\left( z \right)}{{{\left( t-z \right)}^{\alpha-k+1 }}}dz,\,}\,\,\alpha\notin \mathbb{N}_0, \\
 & \frac{d^kg(t)}{dt^k},\,\,\,\,\,\,\,\,\,\,\,\,\,\,\,\,\,\,\,\,\,\,\,\,\,\,\,\,\,\,\,\,\,\,\,\,\,\,\,\,\,\,\,\,\alpha \in\mathbb{N}, \\
& g(t)\,\,\,\,\,\,\,\,\,\,\,\,\,\,\,\,\,\,\,\,\,\,\,\,\,\,\,\,\,\,\,\,\,\,\,\,\,\,\,\,\,\,\,\,\alpha=0 
\end{aligned} \right.
\end{equation*}
is a fractional differential operator of Caputo type \cite{KST06}, where $k=[\alpha]$ is the integer part of $\alpha$.

\medskip
Let $f(x,t)$ be a given function. The direct problem for equation (1) is formulated as follows:

\bigskip
{\bf Direct problem.} To find a solution $u(x,t)$ of equation (1) in $\Omega$, which satisfies
\begin{enumerate}
  \item regularity conditions
  \begin{equation*}
  u, \,\,u_{xx},\,\, \partial_{0t}^{\alpha}u\in C(\Omega),\,\,\int\limits_0^1 \sqrt{x} |u(x,t)|dx<+\infty;
  \end{equation*}
  \item boundary conditions
  \begin{eqnarray}\label{eq2-3}
      \lim\limits_{x\rightarrow 0} x u_x(x,t)=0,\,\, u(1,t)=0, \\
    u(x,0) =0, 0\le x\le 1.
      \end{eqnarray}

\end{enumerate}

Let now $f(x,t)=g(x)$ be an unknown function. In this case, the inverse source problem for equation (1) can be formulated as follows:

\bigskip
{\bf Inverse source problem.} To find a pair of functions $\left\{u(x,t),\, g(x)\right\}$ which
\begin{enumerate}
  \item satisfies regularity conditions
  \begin{equation*}
  u, \,\,u_{xx},\,\, \partial_{0t}^{\alpha}u\in C(\Omega),\,\,\int\limits_0^1 \sqrt{x} |u(x,t)|dx<+\infty;
  \end{equation*}
  \item and in $\Omega$ satisfies equation
  \begin{equation}\label{eq4}
  \partial_{0t}^{\alpha}u(x,t)+\lambda\, \partial_{0t}^{\beta}u(x,t)=u_{xx}(x,t)+\frac{1}{x}u_x(x,t)-\frac{\nu^2}{x^2}u(x,t)+g(x),
  \end{equation}
  \item together with boundary conditions
  \begin{eqnarray}\label{eq5-7}
      \lim\limits_{x\rightarrow 0} x u_x(x,t)=0,\,\, u(1,t)=0, \\
    u(x,0) = \psi(x),\,\, 0\le x\le 1, \\
    u(x,T) = \varphi(x),\,\,0\le x\le 1.
  \end{eqnarray}
  \end{enumerate}
Here $\psi(x),\, \varphi(x)$ are given functions such that
$$
|\psi(0)|<+\infty,\, |\varphi(0)|<+\infty,\,\,\,\psi(1)=0,\,\varphi(1)=0.
$$

\section{Main Results}

In this section we first recall some necessary preliminary facts and then formulate the results.

\subsection{Bessel equation and Fourier-Bessel series}

According to the method of separation of variables, we search for solution of (1) at $f(x,t)\equiv 0$ in the form
\begin{equation*}
u(x,t)=X(x)U(t).
\end{equation*}
Substituting this into (1), with respect to $x$-variable we get the following equation
\begin{equation}\label{eq8}
X''(x)+\frac{1}{x}X'(x)-\frac{\nu^2}{x^2}X(x)=-\mu X(x),
\end{equation}
which is called the Bessel equation of order $\nu$ (see \cite{Wats}). Here $\mu\neq 0$ is a constant.
Conditions (2) give us
\begin{equation}\label{eq9}
\lim\limits_{x\rightarrow 0}x X'(x)=0,\,\,\,X(1)=0.
\end{equation}

The solution of problem (8)-(9) has a form
\begin{equation}\label{eq10}
X(x)=AJ_\nu(\sqrt{\mu}x).
\end{equation}
Substituting this solution to the second condition of (9), we obtain condition for the existence of non-trivial solution of (8) as
\begin{equation}\label{eq11}
J_\nu(\sqrt{\mu})=0.
\end{equation}
It is known that at $\nu>0$, the Bessel function $J_\nu(z)$ has countable number of zeros, moreover, they are real and have pairwise opposite signs. Hence, equation (11) has a countable number of real roots. Denoting $m$th positive root of the equation $J_\nu(z)=0$ by $\gamma_m$, we get $\mu_k=\gamma_k^2,\,k=1,2,...$, for which nontrivial solutions of the eigenvalue problem exist. For arbitrary large $k$, we have (see \cite{Tols})
$$
\gamma_k=k\pi+\frac{\nu\pi}{2}-\frac{\pi}{4}.
$$

Setting in (10) $\mu=\mu_k$, $A=c_k\ne 0$, we get nontrivial solutions (eigenfunctions)
\begin{equation*}
X_k(x)=c_kJ_\nu(\gamma_k x),
\end{equation*}
corresponding to eigenvalues $\gamma_k^2$.

Below we recall some statements from \cite{Tols}:

\textbf{Statement 1.} Functions $J_\nu(\gamma_k x)$ are linearly independent;

\textbf{Statement 2.} For any $l>0$ the equality
$$
\int\limits_0^l x J_\nu \left(\gamma_k \frac{x}{l}\right) J_\nu \left(\gamma_m \frac{x}{l}\right) dx=
\left\{
\begin{array}{cc}
  0 & if \,\,\,\,\,m\ne k;\\
  l^2J_{\nu+1}^2(\gamma_m)/2 & if\,\,\,\,\,\, m=k
\end{array}
\right.
$$
holds true.

\textbf{Statement 3.} System of functions $\left\{\sqrt{x}J_\nu\left(\gamma_k x/l\right)\right\}_{k=1}^\infty$ is complete in $L^2(0,1)$.

\textbf{Statement 4.} The series
\begin{equation}\label{eq12}
\sum\limits_{k=1}^\infty a_k J_\nu\left(\gamma_k \frac{x}{l}\right),
\end{equation}
whose coefficients are defined as
\begin{equation*}
a_k=\frac{2}{l^2J_{\nu+1}^2(\gamma_k)}\int\limits_0^l x\, f(x)J_\nu\left(\gamma_k \frac{x}{l}\right)dx, \,\,k=1,2,...
\end{equation*}
is called the Fourier-Bessel series for a function $f(x)$ on $[0,l]$. If a function is piecewise continuous on  $[0,l]$ and $\int\limits_0^l\sqrt{x}|f(x)|dx<+\infty$, then for $\nu>-1/2$, series (12) converges in every point $x_0\in (0,l)$.

According to these statements, eigenfunctions $X_k(x)=J_\nu(\gamma_k x),\,k=1,2,...,$ of problem (8)-(9) are linearly independent, pairwise orthogonal and any continuous on $(0,1)$ function $f(x)$, satisfying condition $\int\limits_0^1\sqrt{x}|f(x)|dx<+\infty$, can be expanded in $[0,1]$ by the Fourier-Bessel series (see \cite{Tols}).

\subsection{General solution of fractional ODE}

According to the statement of the previous subsection, a solution of the direct problem will be represented as follows:
\begin{equation}\label{eq13}
u(x,t)=\sum\limits_{n=1}^\infty U_n(t) J_\nu\left(\gamma_n x\right),
\end{equation}
\begin{equation}\label{eq14}
f(x,t)=\sum\limits_{n=1}^\infty f_n(t) J_\nu\left(\gamma_n x\right),
\end{equation}
where
\begin{equation}\label{eq15}
U_n(t)=\frac{2}{J_{\nu+1}^2(\gamma_n)}\int\limits_0^1 u(x,t)\, xJ_\nu\left(\gamma_n x\right)dx,
\end{equation}
\begin{equation}\label{eq16}
f_n(t)=\frac{2}{J_{\nu+1}^2(\gamma_n)}\int\limits_0^1 f(x,t)\, xJ_\nu\left(\gamma_n x\right)dx.
\end{equation}

Substituting (13)-(14) into equation (1), we obtain
\begin{equation}\label{eq17}
\partial_{0t}^{\alpha}U_n(t)+\lambda\, \partial_{0t}^{\beta}U_n(t)+\gamma_n^2 U_n(t)=f_n(t).
\end{equation}
General solution of (17) can be written (\cite{KST06}) as
\begin{equation}\label{eq18}
  U_n(t)=\int\limits_0^t f_n(z)(t-z)^{\alpha-1}G_{\alpha,\beta;\lambda,\gamma_n}(t-z)dz+C_1U_{n1}(t),
\end{equation}
where
$$
U_{n1}(t)=\sum\limits_{n=0}^\infty \frac{\gamma_n^{2n}}{n!}t^{\alpha n}{}_1\Psi_1\left[\left.
                                                                                  \begin{array}{l}
                                                                                    (n+1,1) \\
                                                                                    (\alpha n+1, \, \alpha-\beta \\
                                                                                  \end{array}
\,\,\right |\,\,-\lambda t^{\alpha-\beta}
\right],
$$
$$
  G_{\alpha,\beta;\lambda,\gamma_n}(z)=\sum\limits_{n=0}^\infty \frac{\gamma_n^{2n}}{n!}z^{\alpha n}{}_1\Psi_1\left[\left.
                                                                                  \begin{array}{l}
                                                                                    (n+1,1) \\
                                                                                    (\alpha n+\alpha,\, \alpha-\beta \\
                                                                                  \end{array}
\,\,\right |\,\,-\lambda z^{\alpha-\beta}
\right],
$$
$$
{}_1\Psi_1\left[\left.
                                                                                  \begin{array}{l}
                                                                                    (n+1,1) \\
                                                                                    (\alpha n+\beta,\, \alpha) \\
                                                                                  \end{array}
\,\,\right |\,\,z
\right]=\sum\limits_{j=0}^\infty \frac{\Gamma(n+j+1)}{\Gamma(\alpha n+\beta+\alpha j)}\frac{z^j}{j!}=\left(\frac{\partial}{\partial z}\right)^n E_{\alpha,\beta}(z),
$$
and where $E_{\alpha,\beta}(z)$ is two parameter Mittag-Leffler function \cite{KST06}
$$
E_{\alpha,\beta}(z)=\sum\limits_{n=0}^\infty \frac{z^n}{\Gamma(\alpha n+\beta)},\,\,\Re(\alpha)>0,\,\Re(\beta)>0,\,z\in \mathbf{C}.
$$

This representation can be written in a more compact  (and more suitable for further usages) form using the Srivastava-Daoust series in two variables (\cite{sridao}) so that solution (18) can be rewritten as
$$
\begin{array}{l}
U_n(t)=\int\limits_0^t f_n(z) (t-z)^{\alpha-1} S_{1:0;0}^{1:0;0}\left(\aligned
                                              & [1:1;1]: & -; & -; \\
                                             & [\alpha:\alpha;\alpha-\beta]: & -; & -; \\
                                              \endaligned \,|\mu z^{\alpha}, \lambda z^{\alpha-\beta}\right) dz+\\
                                              +C_1 S_{1:0;0}^{1:0;0}\left(\aligned
                                              & [1:1;1]: & -; & -; \\
                                             & [1:\alpha;\alpha-\beta]: & -; & -; \\
                                              \endaligned \,|\mu t^{\alpha}, \lambda t^{\alpha-\beta}\right),
\end{array}
$$
where, the Srivastava-Daoust generalisation of the Kamp\'{e} de F\'{e}riet hypergeometric function in two variables is defined by (\cite[p. 199]{sridao})
$$
\aligned
& S_{C:D;D'}^{A:B;B'} \left(
                                                  \begin{array}{c}
                                                    x \\
                                                    y \\
                                                  \end{array}
                                                \right)
             =S_{C:D;D'}^{A:B;B'}
                \left(\aligned
            &\,\, [(a):\theta,\phi]:\\
            &[(c):\mu,\epsilon]:
            \endaligned\right.\\
              & \hskip 45mm   \left.
            \aligned
            &[(b):\psi];[(b'):\psi']; \\
            &[(d):\eta];[(d'):\eta'];
            \endaligned \, x,y\right)\\
        &  \hskip 35mm   =\sum_{m,n=0}^{\infty}\, \Omega (m,n)\,
                \frac{x^{m}}{m!}  \frac{y^{n}}{n!},
        \endaligned
     $$
where, for convenience,
$$
 \Omega(m,n)=\frac{\displaystyle{\prod_{j=1}^{A}}\Gamma(a_{j}+m\,\theta_{j}+n\,\phi_{j})
            \,\prod_{j=1}^{B}\Gamma(b_{j}+m\,\psi_{j})\,\prod_{j=1}^{B'}\Gamma(b'_{j}+n\,\psi'_{j})}
            {\displaystyle\displaystyle{\prod_{j=1}^{C}}\Gamma(c_{j}+m\,\mu_{j}+n\,\epsilon_{j})
            \,\prod_{j=1}^{D}\Gamma(d_{j}+m\,\eta_{j})\,\prod_{j=1}^{D'}\Gamma(d'_{j}+n\,\eta'_{j})},
$$
and the parameters satisfy
   $$\theta_{1},\cdots,\theta_{A}, \eta'_{1}, \cdots,\eta'_{D}>0. $$

   For convenience, here $(a)$ denotes the sequence of $A$ parameters $a_{1},\cdots, a_{A}$ with similar interpretations for $b_{1}, \cdots, d'$. Empty products should be interpreted as
equal to one. Srivastava and Daoust reported that the series in that representation converges
absolutely for all $x, y\in\mathbb{C}$ when
$$
 \Delta_{1}\equiv 1+\sum_{j=1}^{C}\mu_{j}+\sum_{j=1}^{D}\eta_{j}-\sum_{j=1}^{A}\theta_{j}
-\sum_{j=1}^{B}\psi_{j}>0,
$$
$$
 \Delta_{2}\equiv 1+\sum_{j=1}^{C}\epsilon_{j}+\sum_{j=1}^{D'}\eta'_{j}-\sum_{j=1}^{A}\phi_{j}
-\sum_{j=1}^{B'}\psi'_{j}>0.
$$

From the two representations given above, we can get a new relation between Srivastava-Daoust function in two variables and known Wright function (\cite{wri}).

There is another representation of general solution to (17) (see \cite{gorluch}), which satisfies initial condition $U_n(0)=0$:
\begin{equation}\label{eq19}
U_n(t)=\int\limits_0^tz^{\alpha-1}E_{(\alpha-\beta,\alpha),\alpha}\left(-\lambda z^{\alpha-\beta},-\gamma_n^2 z^\alpha\right)f_n(t-z)dz,
\end{equation}
where
$$
E_{(\alpha-\beta,\alpha),\rho}\left(x,y\right)=\sum\limits_{k=0}^\infty\sum\limits_{i=0}^k \frac{k!}{i! (k-i)!}\cdot\frac{x^iy^{k-i}}{\Gamma(\rho+\alpha k-\beta i)}
$$
is a particluar case of multivariate Mittag-Leffler function (\cite{gorluch}).

Further, we will use this representation since necessary properties of multinomial Mittag-Leffler are available.

\subsection{Formal solution and the convergence}
Substituting (19) into (13), we get a formal solution as
\begin{equation}\label{eq20}
u(x,t)=\sum\limits_{n=1}^\infty \left[\int\limits_0^tz^{\alpha-1}E_{(\alpha-\beta,\alpha),\alpha}\left(-\lambda z^{\alpha-\beta},-\gamma_n^2 z^\alpha\right)f_n(t-z)dz\right]
 J_\nu\left(\gamma_n x\right).
\end{equation}

In order to prove the convergence of this series, we use estimation of the Mittag-Leffler function, obtained in \cite{yamli} (see Lemma 3.2):
$$
|E_{(\alpha-\beta,\alpha),\rho}(x,y)|\leq \frac{C}{1+|x|}.
$$

There is a special approach to prove the convergence of series corresponding to $u_{xx}(x,t)$. Precisely, we have
$$
u_{xx}(x,t)=\sum\limits_{n=0}^\infty U_n(t)\frac{d^2}{dx^2} \left( J_\nu\left(\gamma_n x\right)\right).
$$
Using that $ J_\nu\left(\gamma_n x\right)$ is a solution of the Bessel equation (8), we can rewrite this expression as
$$
u_{xx}(x,t)=\sum\limits_{n=0}^\infty U_n(t)\left[-\gamma_n^2J_\nu(\gamma_n^2 x)+\frac{\nu^2}{x^2}J_\nu(\gamma_n^2 x)-\frac{1}{x}J_\nu'(\gamma_n^2 x)\right].
$$
More regularity conditions for given functions will be needed for the convergence of the following series
$$
\sum\limits_{n=0}^\infty \gamma_n^2 U_n(t) J_\nu(\gamma_n^2 x).
$$
According to the theory of Fourier-Bessel series \cite{Tols} (see p.231), we need to get estimation of the form
$$
\left|\gamma_n^2 U_n(t)\right|\leq \frac{C}{(\gamma_n^2)^{1+\epsilon}},
$$
where $C$ and $\epsilon$ are a positive constants.

If we impose conditions to the given function $f(x,t)$ such that
\begin{itemize}
\item $f(x,t)$ is differentiable four times with respect to $x$;
\item $f(0,t)=f'(0,t)=f''(0,t)=f'''(0,t)=0,\,\,f(1,t)=f'(1,t)=f''(1,t)=0$;
\item $\frac{\partial^4 f(x,t)}{\partial x}$ is bounded,
\end{itemize}
then we will get required estimate in the form of
$$
\left|\gamma_n^2 \int\limits_0^tz^{\alpha-1}E_{(\alpha-\beta,\alpha),\alpha}\left(-\lambda z^{\alpha-\beta},-\gamma_n^2 z^\alpha\right)f_n(t-z)dz\right|\leq \frac{C}{(\gamma_n^2)^{3/2}}.
$$
The convergence of series corresponding to $\partial _{0t}^\alpha u(x,t)$, $u_x(x,t)$ can be done similarly.

Let us now consider the inverse problem.

Substituting (13) and
$$
g(x)=\sum\limits_{n=1}^\infty g_n J_\nu\left(\gamma_n x\right),
$$
where
$$
g_n=\frac{2}{J_{\nu+1}^2(\gamma_n)}\int\limits_0^1 g(x)\, xJ_\nu\left(\gamma_n x\right)dx,
$$
into the equation (4), we obtain
$$
\partial_{0t}^{\alpha}U_n(t)+\lambda\, \partial_{0t}^{\beta}U_n(t)+\gamma_n^2 U_n(t)=g_n.
$$
Solution to this equation, satisfying initial condition $U_n(0)=\psi_n$ has a form
\begin{equation}\label{eq21}
\begin{array}{l}
U_n(t)=g_n \int\limits_0^tz^{\alpha-1}E_{(\alpha-\beta,\alpha),\alpha}\left(-\lambda z^{\alpha-\beta},-\gamma_n^2 z^\alpha\right)dz+\\
+\psi_n\left[1-\lambda t^{\alpha-\beta}E_{(\alpha-\beta,\alpha),1+\alpha-\beta}\left(-\lambda t^{\alpha-\beta},-\gamma_n^2 t^\alpha\right)-\gamma_n^2 t^\alpha E_{(\alpha-\beta,\alpha),1+\alpha}\left(-\lambda t^{\alpha-\beta},-\gamma_n^2 t^\alpha\right)\right],\\
\end{array}
\end{equation}
where
$$
\psi_n=\frac{2}{J_{\nu+1}^2(\gamma_n)}\int\limits_0^1 \psi(x)\, xJ_\nu\left(\gamma_n x\right)dx.
$$
Using the formula
$$
\int\limits_0^tz^{\alpha-1}E_{(\alpha-\beta,\alpha),\alpha}\left(-\lambda z^{\alpha-\beta},-\gamma_n^2 z^\alpha\right)dz=t^\alpha E_{(\alpha-\beta,\alpha),\alpha+1}\left(-\lambda t^{\alpha-\beta},-\gamma_n^2 t^\alpha\right) ,
$$
we rewrite solution (21) as follows
\begin{equation}\label{eq22}
\begin{array}{l}
U_n(t)=\left[g_n-\gamma_n^2\psi_n\right]t^\alpha E_{(\alpha-\beta,\alpha),1+\alpha}\left(-\lambda t^{\alpha-\beta},-\gamma_n^2 t^\alpha\right)+\\
+\psi_n\left[1-\lambda t^{\alpha-\beta}E_{(\alpha-\beta,\alpha),1+\alpha-\beta}\left(-\lambda t^{\alpha-\beta},-\gamma_n^2 t^\alpha\right)\right].\\
\end{array}
\end{equation}
Now we use over-determining condition (7), which passes to $U_n(T)=\varphi_n$ and find unknown coefficient $g_n$ from known $\psi_n$ and $\varphi_n$:
$$
\begin{array}{l}
g_n=\gamma_n^2 \psi_n+\frac{1}{T^\alpha E_{(\alpha-\beta,\alpha),1+\alpha}\left(-\lambda T^{\alpha-\beta},-\gamma_n^2 T^\alpha\right)}\left[\varphi_n-\right.\\
\left.-\psi_n\left[1-\lambda t^{\alpha-\beta}E_{(\alpha-\beta,\alpha),1+\alpha-\beta}\left(-\lambda T^{\alpha-\beta},-\gamma_n^2 T^\alpha\right)\right]\right].\\
\end{array}
$$
Substituting (22) into (13), we get a formal solution as
\begin{equation}\label{eq23}
\begin{array}{l}
u(x,t)=\sum\limits_{n=0}^\infty J_\nu\left(\gamma_n x\right)\cdot \left\{\frac{t^\alpha E_{(\alpha-\beta,\alpha),1+\alpha}\left(-\lambda t^{\alpha-\beta},-\gamma_n^2 t^\alpha\right)}{T^\alpha E_{(\alpha-\beta,\alpha),1+\alpha}\left(-\lambda T^{\alpha-\beta},-\gamma_n^2 T^\alpha\right)}\varphi_n+\right.\\
\left.+\psi_n\left[1-\lambda t^{\alpha-\beta}E_{(\alpha-\beta,\alpha),1+\alpha-\beta}\left(-\lambda t^{\alpha-\beta},-\gamma_n^2 t^\alpha\right)\right]\left[1-\frac{t^\alpha E_{(\alpha-\beta,\alpha),1+\alpha}\left(-\lambda t^{\alpha-\beta},-\gamma_n^2 t^\alpha\right)}{T^\alpha E_{(\alpha-\beta,\alpha),1+\alpha}\left(-\lambda T^{\alpha-\beta},-\gamma_n^2 T^\alpha\right)}\right]\right\}.\\
g(x)=\sum\limits_{n=0}^\infty J_\nu\left(\gamma_n x\right)\cdot\left\{\gamma_n^2 \psi_n+\frac{1}{T^\alpha E_{(\alpha-\beta,\alpha),1+\alpha}\left(-\lambda T^{\alpha-\beta},-\gamma_n^2 T^\alpha\right)}\times\right.\\
\left.\left[\varphi_n-\psi_n\left[1-\lambda t^{\alpha-\beta}E_{(\alpha-\beta,\alpha),1+\alpha-\beta}\left(-\lambda T^{\alpha-\beta},-\gamma_n^2 T^\alpha\right)\right]\right]\right\}.\\
\end{array}
\end{equation}

The convergence part can be done similarly to the case of the direct problem. We note that in order to prove the convergence of series corresponding to $\partial_{0t}^\alpha u(x,t)$, one needs the following formula
$$
\frac{d}{dt}\left[t^{\alpha} E_{(\alpha-\beta,\alpha),\alpha+1}\left(-\lambda t^{\alpha-\beta},-\gamma_n^2 t^\alpha\right)\right]=\lambda\gamma_n^2 t^{\alpha-1}E_{(\alpha-\beta,\alpha),\alpha}\left(-\lambda t^{\alpha-\beta},-\gamma_n^2 t^\alpha\right),
$$
which is proved in \cite{yamli} (see Lemma 3.3).

Expression in the denominator in (23) cannot be equal to zero, but this does not restrict value of $T$ too much.

A uniqueness of solution for direct and inverse problems can be done in the standard way based using the completeness and basis property of the appropriate system.

We can formulate now our results as theorems.

\medskip
{\bf Theorem 1.}
{\em If
\begin{itemize}
\item $f(x,t)$ is differentiable four times with respect to $x$;
\item $f(0,t)=f'(0,t)=f''(0,t)=f'''(0,t)=0,\,\,f(1,t)=f'(1,t)=f''(1,t)=0$;
\item $\frac{\partial^4 f(x,t)}{\partial x}$ is bounded;
\item $f(x,t)$ is continuous and continuously differentiable with respect to $t$,
\end{itemize}
then there exist unique solution of problem (1)-(3) and it is represented by (20).}

\bigskip
{\bf Theorem 2.}
{\em If
\begin{itemize}
\item $\psi(x)$ and $\varphi(x)$ are differentiable four times;
\item $\psi(0)=\psi'(0)=\psi''(0)=\psi'''(0)=0,\,\,\psi(1)=\psi'(1)=\psi''(1)=0$;
\item $\varphi(0)=\varphi'(0)=\varphi''(0)=\varphi'''(0)=0,\,\,\varphi(1)=\varphi'(1)=\varphi''(1)=0$;
\item $\psi^{(4)}(x)$ and $\varphi^{(4)}(x)$ are bounded;
\item $
E_{(\alpha-\beta,\alpha),1+\alpha}\left(-\lambda T^{\alpha-\beta},-\gamma_n^2 T^\alpha\right)\neq 0,
$
\end{itemize}
then there exist unique solution of problem (4)-(7) and it is represented by (23).}

\section{Example}

For the direct problem we choose given function as $f(x,t)=t(x-x^2)^4$ and $\nu=1$. One can easily verify that chosen function satisfies all conditions of the Theorem 1.  In this case, $f_n(t)$ will have this form
$$
f_n(t)=\frac{2}{J_{\nu+1}^2(\gamma_n)}\int\limits_0^1 t x (x-x^2)^4 J_\nu(\gamma_n x)dx=...
$$
We consider the case, when $\beta=0.4,\,\alpha=0.8$ and $\beta=0,\alpha=1$ as an integer case in order to compare fractional and integer cases.

\section{Acknowledgement} This work was done during the visit of authors to ICMS in Edinburgh in July 2016 and was supported by ``Research in Group'' activity of ICMS.
The last author was also supported in parts by the EPSRC
 grant EP/K039407/1 and by the Leverhulme Grant RPG-2014-02.

\end{document}